\documentclass[11pt]{amsart}
\usepackage{graphicx}
\usepackage[mathcal]{euscript}
\usepackage{float}

\restylefloat{figure}

\newtheorem{theorem}{Theorem}[section]
\newtheorem{lemma}[theorem]{Lemma}
\newtheorem{corollary}[theorem]{Corollary}
\newtheorem{proposition}[theorem]{Proposition}

\theoremstyle{definition}

\theoremstyle{remark}

\numberwithin{equation}{section}

\newcommand{\bb}{{\boldsymbol{b}}}

\def\wt{\widetilde}
\def\wh{\widehat}
\def\bb{\mathbb}
\def\OS{Ozsv\'ath-Szab\'o}

\begin{document}

\title{Turaev genus, knot signature, and the knot homology concordance invariants}

\author{Oliver T. Dasbach}
\address{Department of Mathematics\\
Louisiana State University\\
Baton Rouge, Louisiana}
\email{kasten@math.lsu.edu}
\thanks{The first author was partially supported by {NSF-DMS} 0806539 and {NSF-DMS FRG} 0456275.\\
The second author was partially supported by {NSF-DMS} 0739382 ({VIGRE}) and {NSF-DMS} 0602242 ({VIGRE})}
\author{Adam M. Lowrance}
\address{Department of Mathematics\\
University of Iowa\\
Iowa City, Iowa} 
\email{alowrance@math.uiowa.edu}

\subjclass{}
\date{}

\begin{abstract}
We give bounds on knot signature, the Ozsv\'ath-Szab\'o $\tau$ invariant, and the Rasmussen $s$ invariant in terms of the Turaev genus of the knot. 
\end{abstract}

\maketitle

\section{Introduction}

Alternating knots have particularly simple reduced Khovanov homology and knot Floer homology. Lee \cite{Lee:KhovanovAlternating} showed that the reduced Khovanov homology of an alternating knot $K$ is fully determined by its Jones polynomial $V_K(q)$ and its signature $\sigma(K)$. 
Analogously, Ozsv\'ath and Szab\'o \cite{OzsvathSzabo:Alternating} proved that the knot Floer homology of an alternating knot $K$ is determined by its Alexander polynomial $\Delta_K(t)$ and its signature $\sigma(K)$. Furthermore, for alternating knots the Ozsv\'ath-Szab\'o $\tau$ invariant \cite{OzsvathSzabo:FourBallGenus} and the Rasmussen $s$ invariant \cite{Rasmussen:KhovanovSlice} coincide and  are easily computable. In particular, if $K$ is an alternating knot, then
$$2 \tau(K) = s(K) = -\sigma(K).$$
Note, that it took some efforts to show that in general $2 \tau (K)$ and $s(K)$ are not equal \cite{HeddenOrding:ConcordanceInvariantsNotEqual}.

To compute the signature, if $D$ is a reduced alternating diagram of a knot $K$, Traczyk \cite{Traczyk:AlternatingSignature} proved that
\begin{eqnarray*}
\sigma(K) & = & s_A(D) - n_+(D) -1\\
& = & 1 + n_-(D) - s_B(D),
\end{eqnarray*}
where $s_A(D)$ and $s_B(D)$ are the number of components in the all $A$ and all $B$ Kauffman resolutions of $D$ respectively, and $n_+(D)$ and $n_-(D)$ are the number of positive and negative crossings in $D$ respectively. Throughout this paper we choose our sign convention for the signature such that the signature of the positive trefoil is $-2$.

Our goal is to generalize those results to non-alternating knots.
We will examine the relationship between Traczyk's combinatorial knot diagram data and each of the knot signature, the Ozsv\'ath-Szab\'o $\tau$ and the Rasmussen $s$ invariant for all knots. These relationships lead to new lower bounds for the Turaev genus of a knot.

For a given knot diagram in the plane, Turaev \cite{Turaev:SimpleProof} constructed an embedded oriented surface $\Sigma_D$ on which the knot projects. In \cite{DFKLS:KauffmanDessins} it is pointed out that the knot projection is alternating on the Turaev surface and that Turaev surface is a Heegaard surface for $S^3$. The precise construction of the Turaev surface is given in Section \ref{sec::turaev}. The Turaev genus of a knot $g_T(K)$ is the minimum genus of $\Sigma_D$ over all diagrams of the knot.  We will relate the Turaev genus of a knot $K$ with $\sigma(K), \tau(K)$ and $s(K)$ in the following:

\begin{theorem}
\label{thm::lowerbound}
Let $K$ be a knot. Then
\begin{eqnarray*}
\left |\tau(K)+\frac{\sigma(K)}{2} \right | & \leq & g_T(K),\\
\frac{|s(K) + \sigma(K)|}{2}  & \leq & g_T(K),~\text{and}\\
\left |\tau(K) - \frac{s(K)}{2} \right | & \leq & g_T(K).
\end{eqnarray*}
\end{theorem}

For alternating knots, i.e.~when $g_T(K)=0$, those inequalities reflect the results of Oszv\'ath, Szab\'o and Rasmussen.

Abe \cite{Abe:AlternationNumber}, using work of Livingston \cite{Livingston:OSKnotConcordanceInvariant}, has shown that the three quantities on the left in Theorem \ref{thm::lowerbound} are also lower bounds for the alternation number of a knot, which is the minimum Gordian distance between a given knot and any alternating knot. Examining how the Turaev genus of a knot compares to its alternation number remains an interesting open problem.

The paper is organized as follows. In Section \ref{sec::concord}, we review the constructions of the Ozsv\'ath-Szab\'o $\tau$ invariant and the Rasmussen $s$ invariant. In Section \ref{sec::span}, we show a relationship between the spanning tree complexes for reduced Khovanov homology and knot Floer homology. Section \ref{sec::turaev} is a review of the construction of the Turaev surface and its relationship to the spanning tree complexes. Finally, we show how knot signature fits into the picture in Section \ref{sec::signature}. In Section \ref{sec::3braid}, we compute the bounds of Theorem \ref{thm::lowerbound} for knots obtained as the closure of $3$-braids.

The authors would like to thank Josh Greene for helpful conversations.
 
\section{Knot homology concordance invariants}
\label{sec::concord}
In this section, we recall the definitions of the Ozsv\'ath-Szab\'o $\tau$ invariant \cite{OzsvathSzabo:FourBallGenus} and the Rasmussen $s$ invariant \cite{Rasmussen:KhovanovSlice}.

\subsection{Ozsv\'ath-Szab\'o $\tau$ invariant}

Heegaard Floer homology is an invariant for closed $3$-manifolds defined by Ozsv\'ath and Szab\'o in \cite{Ozsvath-Szabo:HolDisks1} and \cite{Ozsvath-Szabo:HolDisks2}. The Heegaard Floer package gives rise to a concordance invariant, called the \OS~$\tau$ invariant, whose construction is given below.

Suppose $(\Sigma, \alpha, \beta, w, z)$ is a Heegaard diagram subordinate to the knot $K$ in $S^3$. This means  $\Sigma$ is a genus $g$ surface and both 
$\alpha = \{\alpha_1,\dots,\alpha_g\}$ and $\beta=\{\beta_1,\dots,\beta_g\}$ are $g$-tuples of homologically linearly independent, pairwise disjoint, simple closed curves in $\Sigma$. Also, $w$ and $z$ are points in the complement of the $\alpha$ and $\beta$ curves in $\Sigma$ lying in a neighborhood of the curve $\beta_1$ and situated on opposite sides of $\beta_1$. The two sets of curves $\alpha$ and $\beta$ are boundaries of attaching disks and specify handlebodies $U_\alpha$ and 
$U_\beta$ both with boundary $\Sigma$ and $U_\alpha\cup_\Sigma U_\beta \cong S^3$. The knot $K$ can be isotoped onto $\Sigma$ such that it is disjoint from $\beta_2,\dots,\beta_g$, an arc of $K$ runs from the basepoint $w$ to the basepoint $z$, and this arc intersects $\beta_1$ once transversely.

Denote the $g$-fold symmetric product of $\Sigma$ by $\text{Sym}^g(\Sigma)$ and consider the two embedded tori $\bb{T}_\alpha = \alpha_1 \times \cdots \times \alpha_g$ and $\bb{T}_\beta = \beta_1\times \cdots \times \beta_g$. Let $\wh{CF}(S^3)$ denote the $\bb{Z}$-module generated by the intersection points of $\bb{T}_\alpha$ and $\bb{T}_\beta$. The complex $\wh{CF}(S^3)$ can be endowed with a differential that counts pseudo-holomorphics disks in Sym$^g(\Sigma)$ between intersection points of $\bb{T}_\alpha$  and $\bb{T}_\beta$. The homology of $\wh{CF}(S^3)$ is denoted $\wh{HF}(S^3)$ and is isomorphic to $\bb{Z}$ (appearing in homological grading zero).

Ozsv\'ath and Szab\'o \cite{OzsvathSzabo:HolomorphicDisks} and independently Rasmussen \cite{Rasmussen:Floer} proved that a knot $K$ induces a filtration on the chain complex $\wh{CF}(S^3)$. Define $\mathcal{F}(K,m)\subset\wh{CF}(S^3)$ to be the subcomplex generated by intersection points with filtration level less than or equal to $m$. There is an induced sequence of maps
$$\imath_K^m:H_*(\mathcal{F}(K,m))\to H_*(\wh{CF}(S^3))=\wh{HF}(S^3)\cong\bb{Z},$$
that are isomorphisms for all sufficiently large integers $m$. The Ozsv\'ath-Szab\'o $\tau$ invariant is defined as
$$\tau(K) = \min\{m\in\bb{Z}~|~\imath_K^m~\text{is non-trivial}\}.$$
By construction $\tau(K)$ is a knot invariant, and Ozsva\'th and Szab\'o \cite{OzsvathSzabo:FourBallGenus} showed that $\tau(K)$ depends only on the concordance class of $K$.

Also, recall that one can use the filtration $\mathcal{F}(K,m)$ to define the knot Floer homology of $K$, denoted $\widehat{HFK}(K)$, as follows. Define
$$\widehat{HFK}(K) = \bigoplus_{m\in\bb{Z}} H_*(\mathcal{F}(K,m) / \mathcal{F}(K,m-1)).$$
Thus $\widehat{HFK}(K)$ is the homology of the complex $\wh{CFK}(K)$, where $\wh{CFK}(K)$ is generated by intersection points of $\bb{T}_\alpha$ and $\bb{T}_\beta$, but unlike in $\wh{CF}(S^3)$, the differential in $\wh{CFK}(K)$ must preserve filtration level.

\subsection{The Rasmussen $s$ invariant}
Khovanov homology \cite{Khovanov:Homology} is a knot invariant that categorifies the Jones polynomial. Rasmussen
\cite{Rasmussen:KhovanovSlice} used Lee's deformation of Khovanov homology \cite{Lee:Endomorphism} to define a concordance invariant, known as the Rasmussen $s$ invariant, whose construction is described below.

\begin{figure}[h]
\includegraphics[scale=.5]{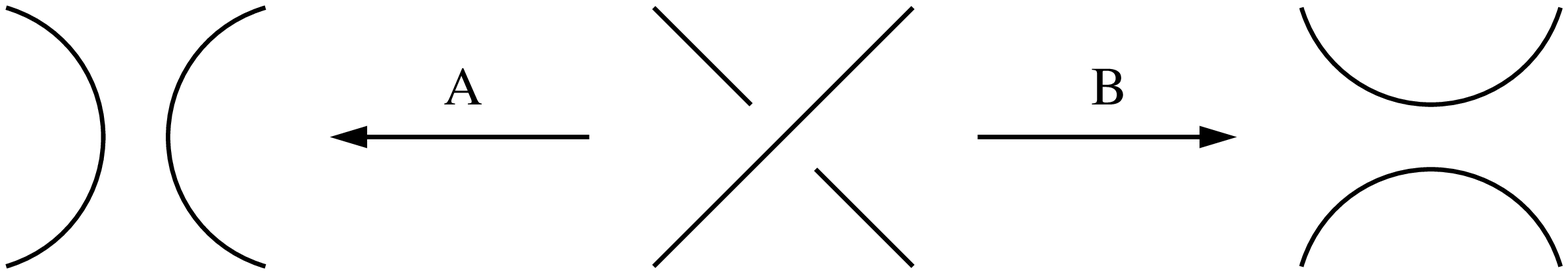}
\caption{The $A$ and $B$ smoothings of a crossing.}
\label{fig::smooth}
\end{figure}

\begin{figure}[h]
\includegraphics[scale=.5]{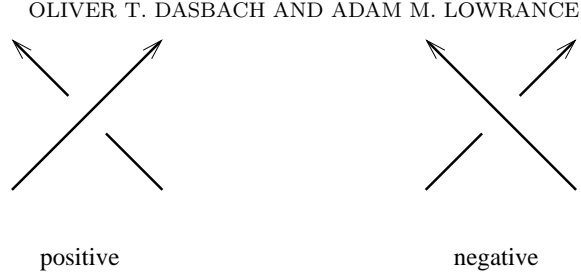}
\caption{The crossing on the left is positive, and the crossing on the right is negative.}
\label{fig::crossing_sign}
\end{figure}

Let $D$ be a diagram of a knot $K$ with crossings labelled $1$ through $k$. Each crossing of $D$ has an $A$-smoothing and a $B$-smoothing, as shown in Figure \ref{fig::smooth}. Associate to each vertex $I$ of the cube $\{A,B\}^k$ the collection of simple closed curves in the plane $D_I$ obtained by smoothing the $i$-th crossing of $D$ according to the $i$-th coordinate of $v$. To each $D_I$ associated the $\bb{Q}$-vector space $V^{\otimes|I|}$ where $V$ is free on two generators $v_+$ and $v_-$ and $|I|$ is the number of components in $D_I$. Define a bigraded $\mathbb{Q}$-vector space, known as the cube of resolutions, by
$$CKh(D)=\bigoplus_{I\in\{A,B\}^k}V^{\otimes|I|}.$$
The homological grading of each summand $V^{\otimes|I|}$ is the number of $B$-smoothings in $I$ minus the number of negative crossings in $D$ (as in Figure \ref{fig::crossing_sign}).

We will investigate two different differentials on $CKh(D)$. The first $\partial_{Kh}$ is Khovanov's differential. The homology $H_*(CKh(D),\partial_{Kh})$ is denoted $Kh(K)$. The vector space $Kh(K)$ has a homological and Jones grading, and its  filtered Euler characteristic is $(q^{1/2}+q^{-1/2})V_K(q)$ where $V_K(q)$ is the Jones polynomial of $K$. (Note that we normalize the Jones grading to be half the usual grading). The second $\partial_{\text{Lee}}$ is Lee's differential. The homology $H_*(CKh(D),\partial_{\text{Lee}})$ is isomorphic to $\bb{Q}\oplus\bb{Q}$. Lee's differential can be written as $\partial_{\text{Lee}}=\partial_{Kh} + \Phi$ where $\Phi$ increases Jones grading. The following theorem is implicit in Lee \cite{Lee:Endomorphism} and explicitly stated in Rasmussen \cite{Rasmussen:KhovanovSlice}.
\begin{theorem}[Rasmussen \cite{Rasmussen:KhovanovSlice}]
\label{thm::spec_seq}
Let $K$ be a knot. There is a spectral sequence with $E_2$ term $Kh(K)$ that converges to $\bb{Q}\oplus\bb{Q}$. 
\end{theorem}

Lee identifies elements of $CKh(D)$ that represent the homology classes $\bb{Q}\oplus\bb{Q}$. These cycles are elements of $V^{\otimes|I|}$ where $I$ is the vertex obtained by smoothing each crossing according the orientation of the knot, i.e.~if a crossing is positive, then one chooses the $A$-smoothing and if a crossing is negative, then one chooses the $B$-smoothing. Therefore, the homological gradings of both of these cycles must be zero. 

Lee's differential does not preserve the Jones grading. In order to obtain a well-defined Jones grading on Lee's homology, one must minimize over all elements in a given homology class. More specifically, if $\alpha\in H_*(CKh(D),\partial_{\text{Lee}})$, then the Jones grading of $\alpha$ is the minimum Jones grading of any element $a$ of $CKh(D)$ such that $a$ represents the homology class $\alpha$.

In \cite{Rasmussen:KhovanovSlice}, Rasmussen showed that Lee's homology is supported in two Jones gradings $s_{\min}(K)$ and $s_{\max}(K)$ depending only on $K$, and moreover $s_{\max} (K)= s_{\min}(K) + 1$. Since our Jones grading is half of Khovanov's original Jones grading, both $s_{\min}(K)$ and $s_{\max}(K)$ are in $\bb{Z}+\frac{1}{2}$. The Rasmussen $s$ invariant is defined as
$$s(K) = s_{\min}(K) + s_{\max}(K).$$
Of course, $s(K)$ is an even integer, and Rasmussen showed that $s(K)$ depends only on the concordance class of $K$.

\section{Spanning tree complexes}
\label{sec::span}

\subsection{Construction of Tait's checkerboard graph}
\begin{figure}[h]
\includegraphics[scale=.4]{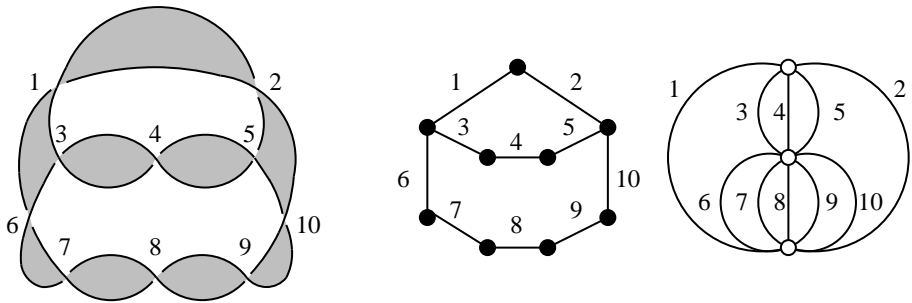}
\caption{A diagram of the $10_{124}$ knot, along with its two Tait graphs. In the black graph, edges $1$ and $2$ are $A$-edges, while edges $3$ through $10$ are $B$-edges. Conversely, in the white graph, edges $1$ and $2$ are $B$-edges, while edges $3$ through $10$ are $A$-edges.}
\label{fig::tait}
\end{figure}

Let $D$ be a diagram of a knot $K$. Color regions of $D$ white and black in a checkerboard fashion, i.e.~so that if two regions are separated by an arc of $D$, then they are different colors. The checkerboard coloring gives rise to the two Tait checkerboard graphs $G$ and $G^*$ of $D$. The vertices of $G$ are in one-to-one correspondence with the black regions, and the edges of $G$ are in one-to-one correspondence with the crossings of $D$. Each edge in $G$ is incident to the vertices that correspond to the black regions near the crossing. An edge in $G$ is called an {\em $A$-edge} (respectively a {\em $B$-edge}) if the $A$-smoothing (respectively the $B$-smoothing) separates the black regions. The vertices of $G^*$ are in one-to-one correspondence with the white regions, and the edges of $G^*$ are in one-to-one correspondence with the crossings of $D$. Each edge in $G^*$ is incident to the vertices that correspond to the white regions near the crossing. If an edge in $G$ is an $A$-edge (respectively a $B$-edge), then the edge corresponding to the same crossing in $G^*$ is a $B$-edge (respectively an $A$-edge). Observe that $G^*$ is the planar dual of $G$. We choose the checkerboard coloring so that the number of $B$-edges in $G$ is greater than or equal to the number of $B$-edges in $G^*$. Figure \ref{fig::tait} shows an example of the Tait graphs for the $10_{124}$ knot. Let $\mathcal{T}(G)$ denote the set of spanning trees of $G$. 

For any subgraph $H$ of $G$, let $V(H)$ be the number of vertices in $H$. Each edge in $G$ is associated to a crossing of $D$, and each crossing in $D$ is either positive or negative (see Figure \ref{fig::crossing_sign}). Moreover, each edge in $G$ is either an $A$-edge or a $B$-edge. For any subgraph $H$ of $G$ or $G^*$, let $E_A^+(H)$ denote the number of edges in $H$ that are both $A$-edges and associated to a positive crossing. Similarly define $E_A^-(H)$, $E_B^+(H)$, and $E_B^-(H)$. 
Also, let $E^+(H)$ denote the number of edges in $H$ associated to positive crossings in $D$ and $E^-(H)$ denote the number of edges in $H$ associated to negative crossings in $D$. Note that $E^+(D) = n_+(D)$ and $E^-(G) = n_-(D)$. Since many of the subsequent arguments rely on graph theoretic ideas, we favor using $E^\pm(G)$ over $n_\pm(D)$. Similarly, let $E_A(H)$ be the number of $A$-edges in $H$ and $E_B(H)$ be the number of $B$-edges in $H$. We alert the reader that in the literature $A$-edges are sometimes called negative edges and $B$-edges are called positive edges. Since we have a different notion of positive and negative edges, we use the $A$ and $B$ notation instead.

If $M=\bigoplus M_{i,j}$ is a finitely generated, bigraded $\bb{Z}$-module, then define the $\delta$-grading of $M$ by $\delta = j - i$.  

\subsection{The knot Floer homology spanning tree complex}
In \cite{OzsvathSzabo:Alternating}, Ozsv\'ath and Szab\'o showed how to associate a Heegaard diagram $(\Sigma, \alpha, \beta, w, z)$ to a knot diagram $D$ such that the intersection points of the tori $\bb{T}_\alpha$ and $\bb{T}_\beta$ embedded into  $\text{Sym}^g(\Sigma)$ are in one-to-one correspondence with the spanning trees of the Tait graph of $D$. Hence there exists a complex whose homology is knot Floer homology that is generated by the spanning trees of the Tait graph.

\begin{proposition}[Ozsv\'ath, Szab\'o \cite{OzsvathSzabo:Alternating}]
\label{prop::span_hfk}
Let $D$ be a diagram of a knot $K$ and let $G$ be its Tait graph. There exists a complex $\wh{CFK}(D)$ whose generators are in one-to-one correspondence with the spanning trees of $G$ and whose homology is $\wh{HFK}(K)$.
\end{proposition}

Ozsv\'ath and Szab\'o \cite{OzsvathSzabo:Alternating} showed how to calculate the $\delta$-grading of a generator by taking a certain sum over the crossings of the knot diagram. In \cite{Lowrance:WidthTuraevGenus}, the second author interpreted the $\delta$-grading in terms of information about the Tait graph of the knot diagram. The 
$\delta$-grading corresponding to a spanning tree $T$ is 
$$\delta_{\wh{HFK}}(T) = \frac{1}{2}\big(E_B^+(T) + E_A^+(G\setminus T) - E_A^-(T) - E_B^-(G\setminus T)\big).$$

\subsection{The Khovanov homology spanning tree complex}
In the cube of resolutions complex for Khovanov homology $CKh(D)$, one associates a two dimensional vector space to each connected component of a Kauffman state. 
Wehrli \cite{Wehrli:SpanningTrees} and Champanerkar and Kofman \cite{ChampanerkarKofman:SpanningTrees} showed that the cube of resolutions $CKh(D)$ retracts onto a complex where one associates a two dimensional vector space to each partial resolution of the knot diagram $D$ that is a twisted unknot (a partial resolution of $D$ that can be transformed into the trivial diagram of the unknot via Reidemeister one moves). The partial resolutions of $D$ that are twisted unknots are in one-to-one correspondence with the spanning trees of the Tait graph of $D$. Similarly, there is a spanning tree complex for reduced Khovanov homology.

Let $G$ be the Tait graph of a knot diagram $D$, and let $\mathcal{T}(G)$ the set of spanning trees of $G$. 
Define the spanning tree complex for Khovanov homology as
$$C(D) = \bigoplus_{T\in\mathcal{T}(G)} \bb{Z}[T_+,T_-],$$
and define the spanning tree complex for reduced Khovanov homology as
$$\wt{C}(D) = \bigoplus_{T\in\mathcal{T}(G)}\bb{Z}[T].$$ 

\begin{proposition}[Wehrli \cite{Wehrli:SpanningTrees}, Champanerkar-Kofman \cite{ChampanerkarKofman:SpanningTrees}]
\label{prop::span_kh}
Let $D$ be a diagram of a knot $K$.
\begin{enumerate}
\item There exists a spanning tree complex $C(D)$ whose homology is $Kh(K)$.
\item There exists a spanning tree complex $\wt{C}(D)$  whose homology is $\wt{Kh}(K)$.
\end{enumerate}
\end{proposition}

Champanerkar and Kofman chose their gradings so that the bigraded Euler characteristic of $\wt{Kh}(K)$ is $q^{-1} V_K(q^2)$ where $V_K(q)$ is the Jones polynomial of $K$. We replace their $j$-grading by $\frac{j+1}{2}$ so that the bigraded Euler characteristic is $V_K(q)$. The gradings between the Khovanov complex and the reduced Khovanov complex are related by
$$i_{Kh}(T_+) = i_{\wt{Kh}}(T) = i_{Kh}(T_-)~\text{and}$$
$$j_{Kh}(T_+)-\frac{1}{2} = j_{\wt{Kh}}(T) = j_{Kh}(T_-) +\frac{1}{2},$$
for any tree $T\in\mathcal{T}(G)$.
The $\delta$-grading corresponding to a spanning tree $T$ in $\wt{C}(D)$ is 
$$\delta_{\wt{Kh}}(T) = E_B(T) + \frac{1}{4}\big(E^+(G) - E^-(G) -E_B(G) + E_A(G) - 2(V(G)-1)\big).$$

For our convenience, we give two alternate formulations of $\delta_{\wt{Kh}}(T)$. Since $T$ is a spanning tree $V(G)-1 = E(T) = E_A(T) + E_B(T),$ and thus
\begin{eqnarray*}
2\delta_{\wt{Kh}}(T) & = & 2 E_B(T) +\frac{1}{2}\big(E^+(G) - E^-(G) -E_B(G) + E_A(G) - 2(E_A(T) + E_B(T))\big)\\
& = & E_B(T) - E_A(T) +\frac{1}{2}\big(E^+(G) - E^-(G) -E_B(G) + E_A(G)\big).
\end{eqnarray*}
The number of crossings of $D$ can be counted in two ways: by counting positive and negative crossings in $D$ and by counting $A$-edges and $B$-edges in $G$. Therefore, $E^+(G) + E^-(G) = E_A(G) + E_B(G)$ or said another way $E^+(G) - E_B(G) = E_A(G) - E^-(G)$. This leads to our two new formulations of $\delta_{Kh}(T)$:
\begin{eqnarray}
\label{eqn::delta1}
2\delta_{\wt{Kh}}(T) & = & E_B(T) - E_A(T) + E^+(G) - E_B(G),~\text{and}\\
\label{eqn::delta2}
2\delta_{\wt{Kh}}(T) & = & E_B(T) - E_A(T) - E^-(G) + E_A(G).
\end{eqnarray}

\subsection{The $\delta$-grading}

The $\delta$-grading of a spanning tree when considered in the reduced Khovanov complex is the same as the $\delta$-grading of that spanning tree when considered in the knot Floer complex. We note that this is not true of the either the homological or polynomial (Jones or Alexander) gradings individually.
\begin{proposition}
\label{delta}
Let $G$ be the Tait graph of a knot diagram $D$. If $T$ is a spanning tree of $G$, then $\delta_{\wt{Kh}}(T) = \delta_{\wh{HFK}}(T)$.
\end{proposition}
\begin{proof}
From equation \ref{eqn::delta1}, we have
\begin{eqnarray*}
2\delta_{\wt{Kh}}(T) & = & E_B(T) - E_A(T) + E^+(G) - E_B(G)\\
& = & E_B^+(T) + E_B^-(T) - E_A^+(T) - E_A^-(T) +E_A^+(G) + E_B^+(G) - E_B^+(G) - E_B^-(G)\\
& = & E_B^+(T) + E_B^-(T) - E_A^+(T) - E_A^-(T) + E_A^+(G) - E_B^-(G)\\
& = & E_B^+(T) - E_B^-(G\setminus T) - E_A^-(T) +  E_A^+(G\setminus T)\\
& = & 2\delta_{\wh{HFK}}(T).
\end{eqnarray*}
\end{proof}

For the remainder of the paper, we use the notation $\delta(T)$ to equivalently mean $\delta_{\wt{Kh}}(T)$ or $\delta_{\wh{HFK}}(T)$. Define
$$\delta_{\min}(D)=\min\{\delta(T)~|~T\in\mathcal{T}(G)\}~\text{and}~
\delta_{\max}(D)=\max\{\delta(T) ~|~T\in\mathcal{T}(G)\}.$$

\begin{proposition}
\label{prop::tau_ineq}
Let $D$ be a diagram of a knot $K$. Then $\delta_{\min}(D)\leq\tau(K)\leq\delta_{\max}(D)$.
\end{proposition}
\begin{proof}
Proposition \ref{prop::span_hfk} implies there is a Heegaard diagram subordinate to $K$ where the intersections points of $\bb{T}_\alpha$ and $\bb{T}_\beta$ are in one-to-one correspondence with the spanning trees of the Tait graph $G$. One can use this Heegaard diagram to generate both the complexes $\wh{CF}(S^3)$ and $\wh{CFK}(K)$.  By the definition of $\tau$, there must be some spanning tree $T$ in filtration level $\tau$. Since the generator of $\wh{HF}(S^3)$ is in homological grading $0$, the tree $T$ must also be in homological grading $0$. Therefore, the tree $T$ (viewed as a generator of $\wh{CFK}(K)$) must satisfy $\delta(T) = \tau(K)$.
\end{proof}

\begin{proposition}
\label{prop::s_ineq}
Let $D$ be a diagram of a knot $K$. Then $2\delta_{\min}(D)\leq s(K)\leq2\delta_{\max}(D)$.
\end{proposition}
\begin{proof}
Since $C(D)$ is a deformation retract of $CKh(D)$, there exists a spectral sequence (analogous to the sequence of Theorem \ref{thm::spec_seq}) whose $E_1$ page is $C(D)$, $E_2$ page is $Kh(K)$ and that converges to $\bb{Q}\oplus\bb{Q}$. Therefore, there exists two generators $T_1$ and $T_2$ of $C(D)$ with $i_{Kh}(T_1)=i_{Kh}(T_2) = 0$ and $j_{Kh}(T_1) = s_{\min}(K)$ and $j_{Kh}(T_2) = s_{\max}(K)$. Hence, there exists a spanning tree $T$ such that $\delta_{\wt{Kh}}(T)=s(K)/2$.
\end{proof}

\section{The Turaev surface}
\label{sec::turaev}

The ideas discussed below involve ribbon graphs associated to a knot diagram. These ideas are developed by Dasbach, Futer, Kalfagianni, Lin, and Stoltzfus (cf. \cite{DFKLS:DDD} and \cite{DFKLS:KauffmanDessins}). The construction of the Turaev surface of a knot diagram is due to Turaev \cite{Turaev:SimpleProof}.

\begin{figure}[h]
\includegraphics[scale=.4]{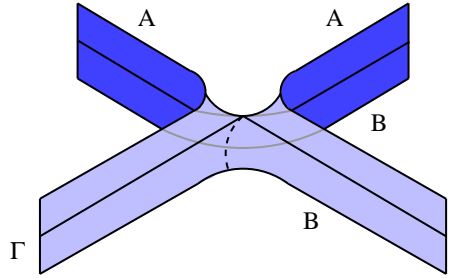}
\caption{In a neighborhood of each vertex of $\Gamma$ a saddle surface transitions between the $A$ and $B$ circles.}
\label{fig::saddle}
\end{figure}

Let $D$ be a knot diagram and $\Gamma$ the $4$-valent plane graph obtained from $D$ by forgetting the ``over-under" information at each crossing. Regard $\Gamma$ as embedded in $\mathbb{R}^2$ which is sitting inside $\mathbb{R}^3$. Remove a neighborhood around each vertex of $\Gamma$, resulting in a collection of arcs in the plane. Replace each arc by a band which is perpendicular to the plane. In the neighborhoods removed earlier, place a saddle so that the circles obtained from choosing an $A$ resolution at each crossing lie above the plane and so that the circles obtained from choosing a $B$ resolution at each crossing lie below the plane. Such a saddle is shown in Figure \ref{fig::saddle}. 
The boundary of the resulting surface is a collection of disjoint circles, where circles corresponding to the all $A$ resolution lie above the plane and circles corresponding to the all $B$ resolution lie below the plane. Cap off each boundary circle with a disk to obtain $\Sigma_D$, the {\em Turaev surface of $D$}. The {\em Turaev genus of a knot $K$} is defined as 
$$g_T(K) = \min\{g(\Sigma_D)~|~D\text{ is a diagram of }K\}.$$

A {\em ribbon graph} is a graph together with a cellular embedding into a surface. The genus $g(\bb{G})$ of a ribbon graph is the genus of the surface into which it embeds.  Denote the number of vertices in a ribbon graph $\bb{G}$ by $V(\bb{G})$. One can embed two ribbon graphs $\bb{A}$ and $\bb{B}$ into $\Sigma_D$ as follows. The vertices of $\bb{A}$ correspond to the disks used to cap off the $A$ circles, and the edges of $\bb{A}$ are the flowlines going from the vertices through the saddles. Similarly, the vertices of $\bb{B}$ correspond to the disks used to cap off the $B$ circles, and the edges of $\bb{B}$ are the flowlines going from the vertices through the saddles. The ribbon graphs $\bb{A}$ and $\bb{B}$ are dual to one another on $\Sigma_D$, and therefore the Euler characteristic of $\Sigma_D$ is determined by
$$\chi(\Sigma_D) = s_A(D) - c(D) + s_B(D),$$
where $c(D)$ is the number of crossings of $D$ and $s_A(D)$ and $s_B(D)$ are the number of components in the all $A$-smoothing and all $B$-smoothing respectively.

Let $\bb{G}$ be a ribbon graph. A {\em ribbon subgraph $\bb{H}$ of $\bb{G}$} is a subgraph of $\bb{G}$ such that the cyclic orientation of the edges in the embedding of $\bb{H}$ is inherited from the embedding of $\bb{G}$. Note that the surfaces on which $\bb{H}$ and $\bb{G}$ are embedded are not necessarily the same. If $\bb{G}$ is embedded on the surface $\Sigma$, then the connected components of $\Sigma\backslash\bb{G}$ are known as the faces of $\bb{G}$. A {\em spanning quasi-tree $\bb{T}$ of $\bb{G}$} is a connected ribbon subgraph of $\bb{G}$ such that $V(\bb{T}) = V(\bb{G})$ and such that $\bb{T}$ has one face. Denote the set of spanning quasi-trees of $\bb{G}$ by $\mathcal{Q}(\bb{G})$.

Recall that $\mathcal{T}(G)$ denotes the set of spanning trees of the Tait graph $G$. Champanerkar, Kofman, and Stoltzfus \cite{CKS:DessinsKhovanov} defined maps $q_{\bb{A}}:\mathcal{T}(G)\to\mathcal{Q}(\bb{A})$ and $q_{\bb{B}}:\mathcal{T}(G)\to\mathcal{Q}(\bb{B})$. Since the sets of edges of $G$, $\bb{A}$, and $\bb{B}$ are each in one-to-one correspondence with the crossings of $D$, we identify all three sets. Because elements of $\mathcal{T}(G)$, $\mathcal{Q}(\bb{A})$, and $\mathcal{Q}(\bb{B})$ are spanning, it suffices to define $q_{\bb{A}}$ and $q_{\bb{B}}$ on the set of edges of $G$. Let $T$ be a spanning tree of $G$.  An $A$-edge of $G$ is in the quasi-tree $q_{\bb{A}}(T)$ if and only if it is in $T$, and a $B$-edge of $G$ is in the quasi-tree $q_{\bb{A}}(T)$ if and only if it is in $G \setminus T$. Similarly, an $A$-edge of $G$ is in the quasi-tree $q_{\bb{B}}(T)$ if and only if it is in $G \setminus T$, and a $B$-edge of $G$ is in $q_{\bb{B}}(T)$ if and only if it is in $T$.
\begin{theorem}[Champanerkar, Kofman, Stoltzfus \cite{CKS:DessinsKhovanov}]
\label{thm::quasi-map}
The maps
$$q_{\bb{A}}: \mathcal{T}(G) \to  \mathcal{Q}(\bb{A})~\text{and}~
q_{\bb{B}}:\mathcal{T}(G)  \to  \mathcal{Q}(\bb{B})$$
are bijections. Moreover, the genera of $q_{\bb{A}}(T)$ and $q_{\bb{B}}(T)$ are determined by
\begin{eqnarray*}
g(q_{\bb{A}}(T)) + E_B(T) & = & \frac{V(G) + E_B(G) -s_A(D)}{2}~\text{and}\\
g(q_{\bb{B}}(T)) + E_A(T) & = & \frac{V(G) + E_A(G) - s_B(D)}{2}
\end{eqnarray*}
\end{theorem}
The following corollary was shown by Champanerkar, Kofman, and Stoltzfus for $\delta_{\wt{Kh}}$ and by the second author for $\delta_{\wh{HFK}}$. In light of Proposition \ref{delta}, it can be seen as a single corollary of the previous theorem. Since the $\delta$-grading for each spanning tree $T$ is the number of $B$-edges in $T$ (up to some overall shift dependent on the diagram $D$), we have the following result.
\begin{corollary}
\label{cor::turaev_genus}
Let $D$ be a knot diagram. The genus of the Turaev surface of $D$ is determined by
$$g(\Sigma_D) = \delta_{\max}(D) - \delta_{\min}(D).$$
\end{corollary}

The maximum and minimum $\delta$-gradings are related to Traczyk's combinatorial data coming from a diagram of the knot.
\begin{corollary}
\label{cor::minmax}
Let $D$ be a knot diagram, and let $G$ be its Tait graph. Then
\begin{eqnarray*}
2\delta_{\min}(D) & =  & s_B(D) - E^-(G) -1~\text{and}\\
2\delta_{\max}(D) & = & 1+E^+(G)-s_A(D).
\end{eqnarray*}
\end{corollary}
\begin{proof}
Let $T_{\min}$ be a spanning tree such that $\delta(T_{\min}) = \delta_{\min}(D)$. By the definition of $q_{\bb{B}}$, the number of edges in $q_{\bb{B}}(T)$ is $E_A(G\setminus T_{\min}) + E_B(T_{\min})$. Since $\delta(T_{\min}) = \delta_{\min}(D)$, the tree $T_{\min}$ has the maximum number of $A$-edges possible, and thus Theorem \ref{thm::quasi-map} implies that $g(q_{\bb{B}}(T_{\min}))=0$. Therefore, $q_{\bb{B}}(T_{\min})$ is a spanning tree of the underlying graph of $\bb{B}$ and has $s_B(D)-1$ edges.

Equation \ref{eqn::delta1} implies
\begin{eqnarray*}
2\delta(T_{\min}) & =  & E_B(T_{\min}) - E_A(T_{\min}) + E_A(G) - E^-(G)\\
& = & E_A(G \setminus T_{\min}) + E_B(T_{\min}) - E^-(G)\\
& = & s_B(D) - E^-(G) - 1. 
\end{eqnarray*}

Similarly, let $T_{\max}$ be a spanning tree such that $\delta(T_{\max})=\delta_{\max}(D)$. By the definition of $q_{\bb{A}}$, the number of edges in $q_{\bb A}(T)$ is $E_A(T) + E_B(G\setminus T)$. Since $\delta(T_{\max})=\delta_{\max}(D)$, the tree $T_{\max}$ has the maximum number of $B$-edges possible, and thus Theorem \ref{thm::quasi-map} implies that $g(q_{\bb{A}}(T_{\max}))=0$. Therefore, $q_{\bb{A}}(T_{\max})$ is a spanning tree of the underlying graph of $\bb{A}$ and has $s_A(D)-1$ edges.

Equation \ref{eqn::delta2} implies
\begin{eqnarray*}
2\delta(T_{\max}) & = & E_B(T_{\max}) - E_A(T_{\max}) +E^+(G) - E_B(G)\\
& = & E^+(G) -E_B(G\setminus T_{\max}) -E_A(T_{\max})\\
& = & 1 + E^+(G) -s_A(D)
\end{eqnarray*}
\end{proof}

\section{Knot signature}
\label{sec::signature}

The signature of a knot $\sigma(K)$ was defined by Trotter in \cite{Trotter:Signature} and was shown to be a concordance invariant by Kauffman and Taylor in \cite{Kauffman:Signature}.
In this section, we show that $\sigma(K)$ satisfies inequalities similar to the inequalities satisfied by $\tau(K)$ and $s(K)$. Consequently, one has new lower bounds for the Turaev genus of a knot.

\subsection{Construction of the Goeritz matrix}

\begin{figure}[h]
\includegraphics[scale=.5]{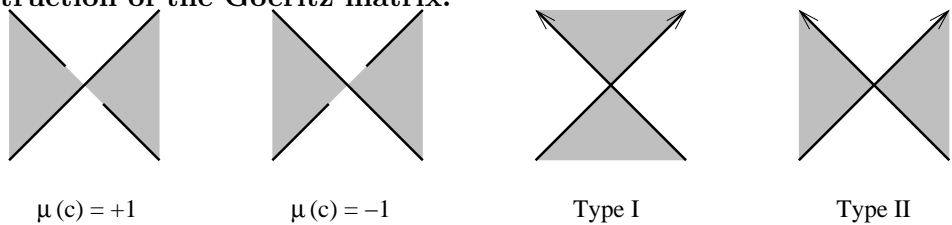}
\caption{The incidence number and type of a crossing.}
\label{fig::goeritz}
\end{figure}

Color the regions of $D$ black and white in a checkerboard fashion. Assume that each crossing is incident to two distinct black regions. Label the black regions of $D$ by $R_0, \dots, R_n$. Assign an incidence number and a type to each crossing, as in Figure \ref{fig::goeritz}. Set 
$$\mu(D) = -\sum_{c~\text{of Type II}}\mu(c).$$
If $i,j\in\{0,1,\dots,n\}$ and $i\neq j$, then define
$$g_{ij} = - \sum_{c\in\overline{R}_i\cap\overline{R}_j}\mu(c),$$
and also, for $i\in\{0,1,\dots,n\}$ define
$$g_{ii}= - \sum_{i\neq j}g_{ij}.$$
Then the Goeritz matrix $\mathcal{G}$ of $D$ is defined to be the $n \times n$ matrix with entries $g_{ij}$ for $i,j\in\{1,\dots,n\}$. Let $\sigma(\mathcal{G})$ denote the signature of the symmetric matrix $\mathcal{G}$, i.e. $\sigma(\mathcal{G})$ is the number of positive eigenvalues $\sigma_+(\mathcal{G})$ minus the number of negative eigenvalues $\sigma_-(\mathcal{G})$. Gordon and Litherland \cite{GordonLitherland:Signature} gave the following formula for the signature of a knot.

\begin{theorem}[Gordon-Litherland \cite{GordonLitherland:Signature}]
Let $D$ be a reduced diagram of a knot $K$. Then $\sigma(K) = \sigma(\mathcal{G}) - \mu(D)$.
\end{theorem}

Observe that the Goeritz matrix is completely determined by the Tait graph $G$. Label the vertices of $G$ by $v_0, v_1, \dots, v_n$ so that the vertex $v_i$ corresponds with the region $R_i$. For $i\neq j$, one can equivalently define
$$g_{ij} = \#(B\text{-edges connecting}~v_i~\text{to}~v_j) - \#(A\text{-edges connecting}~v_i~\text{to}~v_j).$$

\subsection{The $\delta$-grading and signature}

In order to establish the desired inequalities for the signature of the knot, we first need two lemmas.

\begin{lemma}
\label{lem::mirror}
Let $D$ be a knot diagram and $\overline{D}$ its mirror image. Then $\delta_{\min}(D) = -\delta_{\max}(\overline{D})$.
\end{lemma}
\begin{proof}
By Corollary \ref{cor::minmax}, we have
\begin{eqnarray*}
2\delta_{\min}(D) & = & s_B(D) - E^-(G) - 1\\
& = & s_A(\overline{D}) - E^+(\overline{G}) -1\\
& = & - 2\delta_{\max}(\overline{D}).
\end{eqnarray*}
\end{proof}
\begin{lemma}
\label{lem::neg_eig}
Let $D$ be a knot diagram with Tait graph $G$ and Goeritz matrix $\mathcal{G}$. There exists a spanning tree $T\in\mathcal{T}(G)$ such that $E_B(T)\leq\sigma_-(\mathcal{G})$.
\end{lemma}
\begin{proof}
The proof is by induction on the $A$-edges of $G$. First we prove the lemma in the base case where every edge of $G$ is a $B$-edge. Then we show that one can construct the desired spanning tree in $G$ from the graph obtained by contracting an $A$-edge in $G$.

If every edge of $G$ is a $B$-edge, then $D$ is an alternating diagram and the number of components in the all $B$-smoothing $s_B(D)$ is equal to the number of vertices $V(G)$ of $G$. Therefore, the signature of $K$ is given by Traczyk's formula:
$$\sigma(K)=1 + E^-(G)-V(G).$$
Observe that $\mu(D) = E_A^+(G) - E_B^-(G)$,  and hence the Gordon-Litherland formula for signature can be written as
$$\sigma(K) = \sigma_+(\mathcal{G}) - \sigma_-(\mathcal{G}) - E_A^+(G) + E_B^-(G).$$Since each edge of $G$ is a $B$-edge, it follows that $E^-(G) = E_B^-(G)$ and $E_A^+(G)=0$. Therefore,
$$\sigma_+(\mathcal{G}) - \sigma_-(\mathcal{G}) = 1 - V(G).$$
The Goeritz matrix $\mathcal{G}$ is a $(V(G)-1)\times(V(G)-1)$ matrix, and thus $\mathcal{G}$ is negative definite, i.e.~$\sigma_-(\mathcal{G})=V(G)-1$. Hence for any spanning tree $T$ of $G$, we have
$$E_B(T) = V(G) -1 = \sigma_-(\mathcal{G}).$$

Suppose $G$ has $n$ vertices and at least one $A$-edge $e$.  By way of induction, suppose that for all graphs with less than $n$ vertices, there exists a spanning tree $T$ with $E_B(T)$ less than or equal to the number of negative eigenvalues of the Goeritz matrix associated to that graph. Relabel the black regions so that the vertices incident to $e$ are $v_0$ and $v_1$. Let $\mathcal{G}$ be the $n\times n$ Goeritz matrix of $G$ with entries $g_{ij}$, and let $\wt{\mathcal{G}}$ be the $(n-1) \times (n-1)$ Goeritz matrix of the graph $\wt{G}$ obtained by contracting the edge $e$ in $G$ with entries $\wt{g}_{ij}$. Then $\wt{g}_{ij} = g_{i+1,j+1}.$ Therefore $\sigma_-(\wt{\mathcal{G}})\leq\sigma_-(\mathcal{G})$.

By the inductive hypothesis, there exists a spanning tree $\wt{T}$ of $\wt{G}$ such that $E_B(\wt{T})\leq \sigma_-(\mathcal{G})$. One can form a spanning tree $T$ of $G$ by take the edges of $\wt{T}$ and adding the edge $e$. Since $e$ is an $A$-edge, it follows that $E_B(T) = E_B(\wt{T}) \leq \sigma_-(\wt{\mathcal{G}})\leq\sigma_-(\mathcal{G})$. 
\end{proof}

\begin{theorem}
\label{thm::sig_ineq}
Let $D$ be a diagram of a knot $K$. Then $2\delta_{\min}(D)\leq -\sigma(K)\leq2\delta_{\max}(D)$.
\end{theorem}
\begin{proof}
Let $\mathcal{G}$ be the $n \times n$ Goeritz matrix of $D$. By Lemma \ref{lem::neg_eig} there exists is a spanning tree $T\in\mathcal{T}(G)$ such that $E_B(T)\leq\sigma_-(\mathcal{G})$. Since $K$ is a knot, $det(K) = |det(\mathcal{G}| \neq 0$. Therefore $\sigma_-(\mathcal{G}) + \sigma_+(\mathcal{G}) = n = E_A(T) + E_B(T)$. This implies that 
$$0\leq \sigma_-(\mathcal{G}) - E_B(T) = E_A(T) - \sigma_+(\mathcal{G}).$$
Hence $\sigma_+(\mathcal{G})\leq E_A(T)$ and $E_B(T) - E_A(T) \leq \sigma_-(\mathcal{G}) - \sigma_+(\mathcal{G}) = -\sigma(\mathcal{G}).$

Recall that $\mu(D) = E_A^+(G) - E_B^-(G)$. We have
\begin{eqnarray*}
2\delta(T) & = & E_B(T) - E_A(T) + E^+(G) - E_B(G)\\
& = & E_B(T) - E_A(T) + E_A^+(G) + E_B^+(G) - E_B^+(G) - E_B^-(G)\\
& = & E_B(T) - E_A(T) + E_A^+(G) - E_B^-(G)\\
& = & E_B(T) - E_A(T) +\mu(D)\\
&\leq & -\sigma(\mathcal{G}) + \mu(D)\\
& = & -\sigma(K).
\end{eqnarray*}

Therefore, there exists a spanning tree $T$ with $2\delta(T)\leq -\sigma(K)$, and thus for any diagram $D$ of $K$, we have $2\delta_{\min}(D) \leq -\sigma(K)$. 

Let $\overline{D}$ be the mirror image of $D$. By the same argument $2\delta_{\min}(\overline{D})\leq -\sigma(\overline{K})$. By Lemma \ref{lem::mirror}, we have $\delta_{\max}(D) = - \delta_{\min}(\overline{D})$, and of course, $\sigma(K) = -\sigma(\overline{K})$. Therefore $-\sigma(K)\leq 2 \delta_{\max}(D)$.
\end{proof}

\begin{proof}[Proof of Theorem \ref{thm::lowerbound}]
Propositions \ref{prop::tau_ineq} and \ref{prop::s_ineq}, Corollary \ref{cor::minmax}, and Theorem \ref{thm::sig_ineq} imply the following inequalities:
\begin{eqnarray*}
s_B(D) - n_-(D) - 1\leq & 2\tau(K) & \leq 1 + n_+(D) - s_A(D),\\
s_B(D) - n_-(D) - 1\leq & s(K) & \leq 1 + n_+(D) - s_A(D),~\text{and}\\
s_B(D) - n_-(D) - 1\leq & -\sigma(K) & \leq 1 + n_+(D) - s_A(D).
\end{eqnarray*}
The result now follows from Corollary \ref{cor::turaev_genus}.
\end{proof}
 
The third inequality above also follows from Inequality (13.4) in the proof of Theorem 13.3 in \cite{Murasugi:InvariantsOfGraphs} together with results in
\cite{Thistlethwaite:AdequateKauffman}.  

Lobb \cite{Lobb:Rasmussen} gave upper and lower bounds on the Rasmussen $s$ invariant. Lobb's bounds also depend on the diagram of the
knot. He used combinatorial data obtained from the oriented resolution of
the diagram. Our results are similar in nature, but we use combinatorial
data obtain from the all $A$ and all $B$ resolutions.

We conclude this section with a note on unknotting number. Since $|\frac{s(K)}{2}|$, $|\tau(K)|$, and $|\frac{\sigma}{2}|$ are all lower bounds the unknotting number of $K$, the above inequalities give us a way to possibly find a lower bound coming from a diagram of $D$. This lower bound is necessarily weaker than the bounds given by $s(K)$, $\tau(K)$, and $\sigma(K)$.
\begin{proposition}
Let $D$ be the diagram of a knot $K$, and let $G$ be its Tait graph.  Denote the unknotting number of $K$ by $u(K)$.
\begin{enumerate}
\item If $s_B(D) - E^-(G) - 1\geq 0$, then $s_B(D) - E^-(G) -1 \leq 2u(K)$.
\item If $s_A(D) - E^+(G) -1 \geq 0$, then $s_A(D) - E^+(G) -1 \leq 2u(K).$
\end{enumerate}
\end{proposition}

\section{Example: $3$-braid knots}
\label{sec::3braid}

In this section, we examine knots obtained as the closure of a $3$-braid, and compute the bounds of Theorem \ref{thm::lowerbound} for each such knot.

Let $B_3$ denote the braid group on three strands, generated by elements $\sigma_1$ and $\sigma_2$.  Murasugi described the conjugacy classes of closed $3$-braids.
\begin{theorem}[Murasugi \cite{Murasugi:ThreeBraid}]
Any $3$-braid is conjugate to exactly one braid of the form $(\sigma_1\sigma_2)^{3n}\cdot w$, where $n\in\mathbb{Z}$ and $w$ is either
\begin{enumerate}
\item equal to $\sigma_1^{a_1}\sigma_2^{-b_1}\cdots\sigma_1^{a_k}\sigma_2^{-b_k}$, where $a_i,b_i>0$;
\item equal to $\sigma_2^k$ for some $k\in\mathbb{Z}$;
\item equal to $\sigma_1^m\sigma_2^{-1}$ where $m\in\{-1, -2, -3\}.$
\end{enumerate}
\end{theorem}

We say a $3$-braid in one of the above forms is in Murasugi normal form. Closed $3$-braids whose Murasugi normal form is of type (2) or type (3) with $m=-2$ are links. A closed $3$-braid knot of type (3) is a $(3,k)$ torus knot.

\subsection{Torus knots}

Let $T(3,k)$ denote the $(3,k)$ torus knot. Throughout this subsection, we assume $k>0$. The computations for $k<0$ are similar. 
Ozsv\'ath and Szab\'o \cite{OzsvathSzabo:FourBallGenus} and Rasmussen \cite{Rasmussen:KhovanovSlice} computed the value of the $\tau$ and $s$ invariants for torus knots. In our case, we have
$$2\tau(T(3,k))=s(T(3,k))= 2k-2.$$
Gordon, Litherland, and Murasugi  \cite{GLM:Signature} showed that the signature of a $(3,k)$ torus knot is given by
$$\sigma(T(3,6k+l))=-8k-2 l+2,$$
for $l=1,2,4$ or $5$. Therefore, the bounds from Theorem \ref{thm::lowerbound} are
\begin{equation}
\label{eq::torusbound}
\left | \tau(T(3,6k+l))+\frac{\sigma(T(3,6k+l))}{2} \right | =\left| \frac{s(T(3,6k+l))+\sigma(T(3,6k+l))}{2}\right|=2 k,
\end{equation}
where $l=1, 2, 4$ or $5$.

In \cite{Lowrance:TwistedLinks}, the second author found knot diagrams $D_{k,l}$ of the knots $T(3,3k+l)$ such that the genus of the Turaev surface is given by
$$g(\Sigma_{D_{k,l}})=k,$$
where $l=1$ or $2$.
Therefore Equation \ref{eq::torusbound} and Theorem \ref{thm::lowerbound} imply that
$$g_T(T(3,6k+l))=2k,$$
for $l=1$ and $2$. Using other methods, it can be shown that $g_T(T(3,6k+l)) = 2k+1$ for $l=4$ and $5$. In this case, the Equation \ref{eq::torusbound} implies that the bounds from Theorem \ref{thm::lowerbound} are not sharp.

\subsection{Non-torus closed $3$-braids}
We now turn our attention to closed $3$-braid knots whose Murasugi normal form is of type (1). Throughout this subsection, we assume $n>0$. The computations when $n<0$ are similar. Erle  calculated the signature of such a closed $3$-braid knot.
\begin{proposition}[Erle \cite{Erle:Signature}]
If $K_n$ is the closure of $(\sigma_1\sigma_2)^{3n} \sigma_1^{a_1}\sigma_2^{-b_1}\cdots\sigma_1^{a_k}\sigma_2^{-b_k}$, then $$\sigma(K_n) = -4n - \sum_{i=1}^k (a_i - b_i).$$
\end{proposition}

Using work of Van Cott \cite{VanCott:CableKnots}, Greene computed the Rasmussen $s$ invariant for such closed $3$-braids.
\begin{proposition}[Greene \cite{Greene:Unknotting}]
\label{prop::greene}
Let $K_n $ be a knot that is the closure of $(\sigma_1\sigma_2)^{3n} \sigma_1^{a_1}\sigma_2^{-b_1}\cdots\sigma_1^{a_k}\sigma_2^{-b_k}.$ Then
\begin{equation}
s(K_n) = \begin{cases} 6n - 2 -\sigma(K_0), & \text{if $n>0$;}\\
-\sigma(K_0), & \text{if $n=0$;}\\
6n + 2 - \sigma(K_0), & \text{if $n<0.$}
\end{cases}
\end{equation}
\end{proposition}
Greene's proof depends on the following facts.
\begin{enumerate}
\item For a quasi-alternating knot $s(K) = -\sigma(K)$.
\item $s$ is a homomorphism from the smooth knot concordance group $\mathcal{C}\to\bb{Z}$.
\item $|s(K)| \leq 2 g_4(K)$, where $g_4(K)$ is the $4$-genus of $K$.
\item $s$ of the $(m,n)$ torus knot is $(m-1)(n-1)$.
\end{enumerate}
Note that $(2) - (4)$ above are the conditions appearing in Van Cott's \cite{VanCott:CableKnots} work.

Each of $(1) - (4)$ also holds for $2\tau$, and so, using the notation of Proposition \ref{prop::greene}, we have
$$2\tau(K_n) = s(K_n).$$
Therefore
$$\left |\frac{s(K_n) + \sigma(K_n)}{2}\right | =\left |\tau(K_n)+\frac{\sigma(K_n)}{2}\right|=n-1.$$
The second author \cite{Lowrance:TwistedLinks} showed the $g_T(K_n)\leq n$. Hence Theorem \ref{thm::lowerbound} implies
$$g_T(K_n)=n-1~\text{or}~n.$$

\bibliography{../linklit.bib}
\bibliographystyle {amsalpha}

\end{document}